\documentclass[10pt]{article}

\usepackage{amssymb}
\usepackage{amsmath}
\usepackage{theorem}
\usepackage{epsfig}
\usepackage{verbatim}
\usepackage{graphicx}
\usepackage{subfigure}
\usepackage{enumerate}
\usepackage{epstopdf}

\usepackage{color}

\usepackage{indentfirst}
\usepackage{rotating}
\usepackage{listings}
\usepackage{colortbl}
\newtheorem{theorem}{Theorem}[section]

\newtheorem{corollary}[theorem]{Corollary}

\newtheorem{example}[theorem]{Example}

\newcommand{\RR}{\mathbb{R}}

\newcommand{\T}{\mathbb{T}}
\newcommand{\TT}{\mathbb{T}}

\newcommand{\al}{\alpha}

\newcommand{\ep}{\varepsilon}

\newcommand{\pat}{\partial_t}

\newcommand{\paa}{\partial_{\alpha}}

\newcommand{\pa}{\partial}

\newcommand{\be}{\beta}

\newcommand{\intpi}{\int_{-\pi}^\pi}

\newcommand{\alx}{x}

\newenvironment{reduction}[1]{\begin{trivlist} \item[] {\textbf{Reduction of Theorem #1:}}}{\hfill $\Box$
                       \end{trivlist}}

\makeatletter
\renewcommand*\l@section{\@dottedtocline{1}{0em}{1.5em}}
\renewcommand*\l@subsection{\@dottedtocline{2}{1.5em}{2.3em}}
\renewcommand*\l@subsubsection{\@dottedtocline{3}{3.8em}{3.7em}}
\makeatother

\numberwithin{equation}{section}

\allowdisplaybreaks[4]

\begin{document}

\title{Computer-assisted proofs in PDE: a survey}

\author{Javier G\'omez-Serrano}

\maketitle

\begin{abstract}

In this survey we present some recent results concerning computer-assisted proofs in partial differential equations, focusing in those coming from problems in incompressible fluids. Particular emphasis is put on the techniques, as opposed to the results themselves.

\vskip 0.3cm

\textit{Keywords: PDE, computer-assisted, singularity, incompressible}

\end{abstract}

\section{Computer-Assisted Proofs and Interval arithmetics}

In the last 50 years computing power has experienced an enormous development. According to Moore's Law \cite{Moore:moore-law}, every two years the number of transistors has doubled since the 1970's. This phenomenon has resulted in the blooming of new techniques located in the verge between pure mathematics and computational ones. However, even nowadays when we can perform computations at the speeds of the order of Petaflops (a quatrillion floating point operations per second) we can not avoid the following questions, still fundamental in the rigorous analysis of the output of a computer program:

\begin{itemize}
\item[Q1:] Is a computer result influenced by the way the individual operations are done?
\item[Q2:] Does the environment (operating system, computer architecture, compiler, rounding modes, $\ldots$) have any impact on the result?
\end{itemize}

Sadly, the answer to these questions is Yes, which can be easily illustrated by the following C++ codes (see Listings \ref{listing_harmonic} and \ref{listing_rounding}). The first one computes the harmonic series up to a given $N$ in two ways: the first way adds the different numbers from the bigger ones to the smaller and the second one does the sum in the opposite way. The results for $N = 10^{6}$ can be seen in Listing \ref{listing_harmonic_res}. They are not the same and curiously, the real result is not any of the two of them. The second program uses the MPFR library \cite{Fousse-Hanrot-Lefevre-Pelissier-Zimmermann:MPFR} to add two numbers given by the user in two different ways: rounding down and rounding up the result. The output is done in binary. We can see that the results differ (Listing \ref{listing_rounding_res}).

\begin{lstlisting}[label=listing_harmonic,caption=Computation of the truncated Harmonic Series in two different ways]
int main(int argc, char* argv[]){
  int N; cin >> N;
  cout.setf(ios::fixed); cout.precision(15);
  double res1, res2; res1 = res2 = 0.0;
  for (int i=1; i<=N; i++){
    res1 = res1 + 1.0/(double)i;
  }
  for (int i=N; i>=1; i--){
    res2 = res2 + 1.0/(double)i;
  }
  cout << res1 << endl;
  cout << res2 << endl;
}
\end{lstlisting}
\begin{lstlisting}[label=listing_harmonic_res,
caption=Result of the previous computation.]
14.3927267228647811
14.3927267228657563
\end{lstlisting}

\begin{lstlisting}[label=listing_rounding,caption=Sum of two numbers with different rounding]
int main (int argc, char **argv){
  mpfr_t x, y, d, u;
  mpfr_prec_t prec;
  prec = atoi (argv[1]);
  int pprec = prec - 1;
  mpfr_inits2 (prec, x, y, d, u, (mpfr_ptr) 0);
  mpfr_set_str (x, argv[2], 0, GMP_RNDN);
  mpfr_printf ("x = %.*Rb\n", pprec, x);
  mpfr_set_str (y, argv[3], 0, GMP_RNDN);
  mpfr_printf ("y = %.*Rb\n", pprec, y);

  mpfr_add (d, x, y, GMP_RNDD);
  mpfr_printf ("d = %.*Rb\n", pprec, d);

  mpfr_add (u, x, y, GMP_RNDU);
  mpfr_printf ("u = %.*Rb\n", pprec, u);

  return 0;
}
\end{lstlisting}

\begin{lstlisting}[label=listing_rounding_res,caption=Program executed with arguments 10 0.1 1]
x = 1.100110011p-4
y = 1.000000000p+0
d = 1.000110011p+0
u = 1.000110100p+0
\end{lstlisting}

This shows that even the simplest algorithms need a careful analysis: only two operations suffice to give different results if executed in different order or with different rounding methods.

\subsection{What is a Computer-assisted proof?}

A computer-assisted proof may mean different things depending on the field. For simplicity, we will focus in the context of PDE or ODEs. The starting point is typically an object, or in a broader sense, a behaviour. Most computer-assisted proofs are devoted to show an instance of these objects or behaviours. This is usually done via two steps:

\begin{enumerate}
\item[Step 1:] Perform an analytical reduction of the problem into a (possibly big) set of (typically open) conditions in a way that if those conditions are met, then the theorem is proved. Sometimes the amount of conditions is too big or too complicated to be checked by hand, even though a human could do it without computer-assistance given enough time and space. Examples of tour de force human checked conditions can be found in \cite{Costin-Kim-Tanveer:quasi-solution-blasius,Costin-Tanveer:analytical-approximation-blasius}.

\item[Step 2:] Use the computer to rigorously validate the set of conditions from Step 1.
\end{enumerate}

We now present more details of the technique in 3 flavours, outlining some of its possible applications:

\begin{example}\label{example1}
Compute explicit, tight bounds of hard (singular) integrals and use them to track short time behaviour of solutions of a PDE (see \cite{GomezSerrano-GraneroBelinchon:turning-muskat-computer-assisted,Cordoba-GomezSerrano-Zlatos:stability-shifting-muskat-II,Castro-Cordoba-GomezSerrano-MartinZamora:remarks-geometric-properties-sqg}).
\end{example}

\begin{example}\label{example2}
Bound the norm of a given operator, then use a fixed point theorem to show (by contractivity) existence of solutions, even for ill-posed or singular problems (see \cite{DeLaLlave-Sire:a-posteriori-kam-ill-posed-pde,Enciso-GomezSerrano-Vergara:convexity-cusped-whitham,Castelli-Gameiro-Lessard:rigorous-numerics-ill-posed-pde}).
\end{example}

\begin{example}\label{example3}
Track the spectrum of a given operator, and use this information to say something about the stability/instability, solve an eigenvalue problem, or quantify spectral gaps (see \cite{Castro-Cordoba-GomezSerrano:global-smooth-solutions-sqg}).
\end{example}

In this paper we will give instances of examples \ref{example1} and \ref{example3}, in Sections \ref{sectionmuskat} and \ref{sectionsqg} respectively.

A natural concern is whether the intermediate numerical calculations may be error-prone or not, or depend on the implementation. To circumvent this problem, we will give up on calculating the exact answer (since we are now interested in checking inequalities - as opposed to equalities -) and factor in errors, etc.

The theory of interval analysis developed by R. Moore \cite{Moore-Bierbaum:methods-applications-interval-analysis} is an example of a tool, which albeit being impractical due to inefficient resources at the time of its conception, is now being widely used. It belongs to the paradigm known as rigorous computing (in some contexts also called validated computing), in which numerical computations are used to provide rigorous mathematical statements about a result. The philosophy behind the theory of interval analysis consists in working with and producing objects which are not numbers, but intervals in which we are sure that the true result lies. Nevertheless, we should be precise enough since even with plenty of resources, overestimation might lead to too big intervals which might not guarantee the desired result.

In this setting, the main building blocks are intervals containing the answer and intermediate results. In Section \ref{sectionintarith} we explain how to develop an arithmetic for these objects and how to work with it.

Nowadays, there are a few free libraries that implement interval arithmetics and applications such as CAPD \cite{CAPD}, C-XSC \cite{CXSC}, MPFR/MPFI \cite{Fousse-Hanrot-Lefevre-Pelissier-Zimmermann:MPFR,Revol-Rouillier:MPFI} or Arb \cite{Johansson:Arb}.

\subsection{History}

Lately, interval methods have become quite popular among mathematicians. Several highly non-trivial results have been established by the use of interval arithmetics, see for example \cite{Lanford:cap-feigenbaum,Hales:kepler-annals,Hass-Schlafly:double-bubbles-annals,Gabai-Meyerhoff-Thurston:homotopy-hyperbolic-3manifolds,Schwartz:obtuse-triangular-billiards-100-degrees-I} as a small sample. 

In analysis, the most celebrated result is the proof of the dynamics of the Lorenz attractor (Smale's 14th Problem) \cite{Tucker:lorenz-focm}. However, the study of the dynamics of a system has been restricted until very recently almost always to (typically low-dimensional) ODEs. Some examples involving ODEs but an infinite dimensional system are the computation of the ground state energy of atoms or the relativistic stability of matter (see \cite{Seco:thesis-lower-bounds-ground-state-energy-atoms,Fefferman-Seco:aperiodicity-hamiltonial-flow,Fefferman-DeLaLLave:relativistic-stability-of-matter-I}).

Regarding PDE, most of the work has been carried out for dissipative systems (i.e. systems in which the $L^{2}$-norm of the function decreases with time). The most popular ones are the Kuramoto-Sivashinsky equations or Navier-Stokes in low dimensions. The main feature of these models is that one can study the first $N$ modes of the Fourier expansion of the function and see the rest as an ``error''. Since the system is dissipative, if $N$ is large enough, one can get a control on the error throughout time.

Other techniques reduce the problem to compute the norm of an operator (as in Example \ref{example2}) and apply a fixed point theorem, or use topological methods (Conley index). They have been successfully applied for instance in computing the following: Conley index for Kuramoto-Sivashinsky \cite{Zgliczynski-Mischaikow:rigorous-numerics-kuramoto}, bifurcation diagram for stationary solutions of Kuramoto-Sivashinsky \cite{Arioli-Koch:cap-stationary-ks}, stationary solutions of viscous 1D Burgers with boundary conditions \cite{Fogelklou-Tucker-Kreiss-Siklosi:cap-boundary-value-problem-integral-boundary-condition}, traveling wave solutions for 1D Burgers equation \cite{Fogelklou-Tucker-Kreiss:cap-traveling-wave-burgers}, periodic orbits of Kuramoto-Sivashinsky \cite{Figueras-DeLaLLave:cap-periodic-orbits-kuramoto,Zgliczynski:periodic-orbit-kuramoto,Figueras-Gameiro-Lessard-DeLaLLave:framework-cap-invariant-objects,Gameiro-Lessard:periodic-orbits-ks,Arioli-Koch:integration-dissipative-pde}, Conley index for the Swift-Hohenberg equation \cite{Day-Hiraoka-Mischaikow-Ogawa:rigorous-numerics-swift-hohenberg}, bifurcation diagram of the Ohta-Kawasaki equation \cite{vandenBerg-Williams:validation-2d-ohta-kawasaki}, stability of periodic viscous roll waves of the KdV-KS equation \cite{Barker:numerical-proof-stability-roll-waves}, existence of hexagons and rolls for a pattern formation model \cite{vandenBerg-Deschenes-Lessard-Mireles:hexagons-rolls}, self-similar solutions of a 1D model of 3D axisymmetric Euler \cite{Hou-Liu:self-similar-1d-choi-kiselev-yao}, and many others.

In the elliptic setting,  similar techniques have been developed via finite elements \cite{Nagatou-Yamamoto-Nakao:numerical-verification-elliptic-pde,Yamamoto-Nakao:numerical-verifications-elliptic-fem}.
Very recently, there are papers also dealing with the hyperbolic PDE case \cite{Arioli-Koch:periodic-solutions-hamiltonian-pde}.

We also point the reader to the expository article \cite{vandenBerg-Lessard:rigorous-numerics-notices}, to the excellent monographs \cite{Moore-Bierbaum:methods-applications-interval-analysis,Tucker:validated-numerics-book,Rump:verification-methods-survey} and the survey \cite{Nakao:numerical-verification-ode-pde}.

\section{Interval Arithmetics}\label{sectionintarith}

\subsection{Basic Arithmetic}

Representing an abstract concept such as a real number by a finite number of zeros and ones has the advantage that the calculations are finite and the framework is practical. The drawback is naturally that the amount of numbers that can be written in this way is finite (although of the same order of magnitude as the age of the universe in seconds) and inaccuracies might arise while performing mathematical operations. We will now discuss the basics of interval arithmetics. 

Let $\mathbb{F}$ be the set of representable numbers by a computer. We will work with the set of representable closed intervals $\mathbb{IR}
= \{[\underline{a},\overline{a}] | \quad \underline{a} \leq \overline{a}, \quad \underline{a},\overline{a} \in \mathbb{F}\}$. For every element $[a] \in \mathbb{IR}$ we will refer to it by either $[a]$ or by $[\underline{a},\overline{a}]$, whenever we want to stress the importance of the endpoints of the interval. We can now define an arithmetic by the theoretic-set definition
\begin{align}
[x] \star [y] = \{x \star y | \quad x \in [x], y \in [y]\},
\end{align}
for any operation $\star \in \{+,-,\times,\div\}$. We can easily define them by the following equations:
\begin{align*}
[x]+[y] & = [\underline{x}+\underline{y},\overline{x}+\overline{y}] \\
[x]-[y] & = [\underline{x}-\overline{y},\overline{x}-\underline{y}] \\
[x]\times[y] & = [\min\{\underline{x}\underline{y},\underline{x}\overline{y},\overline{x}\underline{y},\overline{x}\overline{y}\},
\max\{\underline{x}\underline{y},\underline{x}\overline{y},\overline{x}\underline{y},\overline{x}\overline{y}\}] \\
[x] \div [y] & = [x] \times \left[\frac{1}{\overline{y}},\frac{1}{\underline{y}}\right], \text{ whenever } 0 \not \in [y].
\end{align*}

Note that this interval-valued operators can be extended to other algebraic expressions involving exponential, trigonometric, inverse trigonometric functions, etc. This derivation is purely theoretical, and we should keep in mind that, if carried out on a computer, the results of an operation have to be rounded up or down according to whether we are calculating the left or right endpoint so that the true result is enclosed in the produced interval. The main feature of the arithmetic is that if $x \in [x], y \in [y]$, then necessarily $x \star y \in [x] \star [y]$ for any operator $\star$. This property is fundamental in order to ensure that the true result is always contained in the interval we get from the computer.

We remark that this arithmetic is not distributive, but subdistributive, i.e:
\begin{align*}
[a]\times([b]+[c]) & \neq [a]\times[b] + [a]\times[c] \\
[a]\times([b]+[c]) & \subset  [a]\times[b] + [a]\times[c]
\end{align*}
\begin{example}
If we set $[a] = [3,4]$, $[b] = [1,2]$, $[c] = [-1,1]$, then:
\begin{align*}
[a]\times([b]+[c]) & = [3,4] \times [0,3] = [0,12]\\
[a]\times[b] + [a]\times[c] & = [3,4] \times [1,2] + [3,4] \times [-1,1] = [3,8] + [-4,4] = [-1,12]
\end{align*}
\end{example}
This illustrates that the way in which operations are executed in the interval-based arithmetic matters much more than in the real-based. As an example, consider the function $f(x) = 1-x^2$ and a domain $D = [-1,1]$. Over the reals, we can write $f$ as any of the following functions:
\begin{align*}
f_1(x) & = 1 - x^2 \\
f_2(x) & = 1 - x \cdot x \\
f_3(x) & = (1+x)\cdot(1-x)
\end{align*}
However, evaluating $f_i$ over $D$ we get the enclosures:
\begin{align*}
f_1(D) & = [0,1] \\
f_2(D) & = [0,2] \\
f_3(D) & = [0,4]
\end{align*}
We observe that although $f_3$ is completely factored, if we expand it we get an expression of the form $x - x$ which in the interval-based arithmetic is equal to an interval of a width twice the width of the domain in which we are evaluating the expression: a price too high to pay compared with the width of the interval $[0,0]$, another form to write the same expression over the reals.

For readability purposes, instead of writing the intervals as, for instance, $[123456,123789]$, we will sometimes instead refer to them as $123^{456}_{789}$.

\subsection{Automatic Differentiation}

One of the main tasks in which we will need the help of a computer is to calculate a massive amount of function evaluations and their derivatives up to a given order at several points and intervals. In order to perform it, one could first think of trying to differentiate the expressions symbolically. However, we don't need the expression of the derivative, just its evaluation at given points. This, together with the fact that the amount of terms of the derivative might grow exponentially with the number of derivatives taken, makes the use of symbolic calculus impractical. Instead of calculating the expression of every derivative, we will use the so-called \emph{automatic differentiation} methods. Suppose $f(x)$ is a sufficiently regular function and let $x_0$ be the point (or interval) of which we want to calculate its image by $f$. We define
\begin{align*}
(f)_{0} & = f(x_{0}) \\
(f)_{k} & = \frac{1}{k!} \frac{d^{k}}{dx^{k}}f(x_0), \quad k = 1,2,\ldots,N,
\end{align*}
where $N$ stands for the maximum number of derivatives of the function we want to evaluate. We can think about $(f)$ as being the coefficients of the Taylor series around $x_0$ up to order $N$. We now show how to compute the coefficients $(f)$ for some of the functions that will appear in our programs. The generalization of the missing functions is immediate. However, it is possible to derive similar formulas for any solution of a differential equation (e.g. Bessel functions).
\begin{align*}
(u\pm v)_{k} & = (u)_{k} \pm (v)_{k} \\
(u \cdot v)_{k} & = \sum_{j=0}^{k}(u)_{j}(v)_{k-j} \\
(u \div v)_{k} & = (1/v) \left((u)_{k} - \sum_{j=1}^{k}(v)_{j}(u \div v)_{k-j}\right) \\
(\sin(u))_{k} & = \frac{1}{k}\sum_{j=0}^{k-1}(j+1)(\cos(u))_{k-1-j}(u)_{j+1} \\
(\cos(u))_{k} & = -\frac{1}{k}\sum_{j=0}^{k-1}(j+1)(\sin(u))_{k-1-j}(u)_{j+1}
\end{align*}

Automatic differentiation has become a natural technique in the field of Dynamical Systems, since the cost for evaluating an expression up to order $k$ is $O(k^2)$, making it a fast and powerful tool to approximate accurately trajectories \cite{Simo:global-dynamics-fast-indicators}. It has also been used for the computation of invariant tori and their associated invariant manifolds \cite{Haro-DeLaLlave:parameterization-invariant-tori-whiskers-explorations,Haro-DeLaLlave:parameterization-invariant-tori-whiskers-rigorous} or the computation of normal forms of KAM tori \cite{Haro:algorithm-canonical-transformations}. For more applications in Dynamical Systems we refer the reader to the book \cite{Haro-Canadell-Figueras-Luque-Mondelo:parameterization-method-book}. Automatic Differentiation is also an important element in the so-called Taylor models \cite{Neumaier:taylor-forms,Makino-Berz:taylor-models,Joldes:rigorous-polynomial-approximations}, in which functions are represented by couples $(P,\Delta)$, being $P$ a polynomial and $\Delta$ an interval bound on the absolute value of the difference between the function and $P$. Nowadays, there are several packages that implement it, for example \cite{Jorba-Zou:taylor-package,Barrio:taylor-series-odes}.

\subsection{Integration}

In this section we will discuss the basics of rigorous integration. A few examples where rigorous integration has been used or developed are \cite{Berz-Makino:high-dimensional-quadrature,Kramer-Wedner:adaptive-gauss-legendre-verified-computation,Lang:multidimensional-verified-gaussian-quadrature}. A more detailed version concerning singular integrals can be found in the next subsection. We will only give the details of the one-dimensional case, omitting the multidimensional one, which can be done extending the methods in a natural way.

The main problem we address here is to calculate bounds for a given integral

\begin{align*}
I = \int_{a}^{b} f(x)dx, \quad -\infty < a < b < \infty.
\end{align*}

Different strategies can be used for this purpose. For instance, we can extend the classical integration schemes:

\begin{align*}
I = \sum_{i=1}^{N} \int_{x_{i-1}}^{x_i}f(x)dx, \quad a = x_0 < x_1 < \ldots < x_{N} = b.
\end{align*}

In every interval, we approximate $f(x)$ by a polynomial $p(x)$ and an error term. We detail some typical examples in Table \ref{tab:rig_integration}:

\renewcommand{\arraystretch}{1.5}
\begin{table}[h]
	\centering
	\begin{tabular}{|c|c|c|c|}
		\hline
		& Midpoint Rule & Trapezoid Rule & Simpson's Rule \\
		\hline
	& \multicolumn{3}{ |c| }{$\displaystyle I \approx \sum_{i=1}^{N}p_i(x)dx$} \\
	\hline
	$p_i(x)$ & $f\left(\frac{x_{i}+x_{i-1}}{2}\right)$ &
	$\begin{array}{l}
	\displaystyle f(x_{i-1})\left(1-\frac{x-x_{i-1}}{h_i}\right) \\
	\displaystyle + f(x_{i})\left(\frac{x-x_{i-1}}{h_i}\right)
	\end{array}$
	  &
	$\begin{array}{l}
	\displaystyle f(x_{i-1})\frac{(x-x_{i})(x-\frac{x_{i}+x_{i-1}}{2})}{h_i^2/2}  \\
	 \displaystyle - f\left(\frac{x_{i}+x_{i-1}}{2}\right)\frac{(x-x_{i})(x-x_{i-1})}{h_i^2/4}  \\
	 \displaystyle + f\left(x_{i}\right)\frac{(x-x_{i-1})(x-\frac{x_{i}+x_{i-1}}{2})}{h_i^2/2}
	\end{array}$ \\
	\hline
	Error & $\frac{b-a}{24}h^2f^{2}([a,b])$ & $-\frac{b-a}{12}h^2f^{2}([a,b])$ & $-\frac{b-a}{2880}h^4f^{4}([a,b])$ \\
	\hline
	\end{tabular}
	\caption{Different Schemes for the rigorous integration.}	
	\label{tab:rig_integration}
\end{table}
\renewcommand{\arraystretch}{1.5}

It is now clear where the interval arithmetic takes place. In order to enclose the value of the integral, we need to compute rigorous bounds for some derivative of the function at the integration region.

Another approach consists of taking the Taylor series of the integrand up to order $n$ as the polynomial $p_i(x)$. Centering the Taylor series in the midpoint of the interval makes us integrate only roughly over half of the terms (since the other half are equal to zero). We can see that

\begin{align*}
\int_{a}^{b}f(x)dx & = \int_{a}^{b}\left(f(a) + (x-a)f'(a)+ \ldots + \frac{(x-a)^{n}}{n!}f^{n}(a)+\frac{(x-a)^{n+1}}{(n+1)!}f^{n+1}(\xi(x))\right)dx \\
& \in \int_{a}^{b}\left(f(a) + (x-a)f'(a)+ \ldots + \frac{(x-a)^{n}}{n!}f^{n}(a)+\frac{(x-a)^{n+1}}{(n+1)!}f^{n+1}([a,b])\right)dx \\
& = \underbrace{(b-a)f(a) + \frac{1}{2}(b-a)^2f'(a) + \ldots +\frac{(x-a)^{n+1}}{(n+1)!}f^{n}(a)}_{\text{Real number (thin interval)}}+\underbrace{\frac{(x-a)^{n+2}}{(n+2)!}f^{n+1}([a,b])}_{\text{Error (thick interval)}}.
\end{align*}

We now compare the two methods in the following examples, in which we integrate $\int_{0}^{1}e^{x}dx$.
\begin{example}
If we take $N = 4$ and use a trapezoidal rule, we enclose the integral in
\begin{align*}
\int_{0}^{1}e^{x}dx & = \frac{1}{2}\left(e^{0} + 2e^{1/4} + 2e^{1/2} + 2e^{3/4} + e^{1}\right)\frac{1}{4} - \frac{1}{12}\frac{1}{16}e^{[0,1]} \\
& \in [1.72722,1.72723] - [0.0050283, 0.014578] = [1.712642,1.7222017]
\end{align*}
\end{example}

\begin{example}
\begin{align*}
\int_{0}^{1} e^{x} dx & \in \int_{0}^{1} 1 + x + \frac{x^{2}}{2} + \frac{x^{3}}{6}e^{[0,1]}dx \\
& = \left.x + \frac{x^{2}}{2} + \frac{x^{3}}{6}\right|^{x=1}_{x=0} + \left.[1,e]\frac{x^{4}}{24}\right|^{x=1}_{x=0} \\
& = \frac{10}{6} + \frac{1}{24}[1,e] \\
& = [1\text{.}70833,1\text{.}77994]
\end{align*}
\end{example}

The exact result is $e - 1 \approx 1.71828182846$. We can see that there is a tradeoff between function evaluations (efficiency of the scheme) and quality (precision) of the results, since the first method is more exact but requires more evaluations of the integrand, while for the second it is enough to compute the Taylor series of the integrand.

\subsection{Singular integrals and integrals over unbounded domains}\label{subsection:singularint}

In this subsection we will discuss the computational details of the rigorous calculation of some singular integrals. In particular we will focus on the Hilbert transform, but the methods apply to any singular integral. Parts of the computation  (the $N$ and $F$ parts) are slightly related to the Taylor models with relative remainder presented in \cite{Joldes:rigorous-polynomial-approximations}. See also the paper \cite{Castro-Cordoba-Fefferman-Gancedo-GomezSerrano:stability-splash-singularities-water-waves} regarding the rigorous inversion of operators involving singular integrals.

Let us suppose that we have a $2\pi$-periodic function $f$, which for simplicity we will assume it is $C^k$ (this requirement can be relaxed). We want to calculate rigorously the Hilbert Transform of $f$, that is

\begin{align*}
Hf(x) = \frac{PV}{\pi} \int_{\T} \frac{f(x)-f(x-y)}{2\tan\left(\frac{y}{2}\right)}dy,
\end{align*}

We can split our integral in

\begin{align*}
Hf(x) & = \frac{PV}{\pi} \int_{\T} \frac{f(x)-f(x-y)}{2\tan\left(\frac{y}{2}\right)}dy \\
& = \frac{PV}{\pi} \int_{|y| < \varepsilon_1} \frac{f(x)-f(x-y)}{2\tan\left(\frac{y}{2}\right)}dy
+ \frac{PV}{\pi} \int_{\varepsilon_1 \leq |y| < \pi - \varepsilon_2} \frac{f(x)-f(x-y)}{2\tan\left(\frac{y}{2}\right)}dy
+ \frac{PV}{\pi} \int_{|y| \geq \pi - \varepsilon_2} \frac{f(x)-f(x-y)}{2\tan\left(\frac{y}{2}\right)}dy \\
& \equiv H^{N} f(x) + H^{C}f(x) + H^{F}f(x).
\end{align*}

The integration of $H^{C}f(x)$ is easy since the integrand is smooth, and the denominator is bounded away from zero.

We now move on the the term $H^{N}f(x)$. In this case, we perform a Taylor expansion in both the denominator
\begin{align*}
2\tan\left(\frac{y}{2}\right) = (y)+c(\varepsilon_1)(y)^{3}, \quad c(\varepsilon_1) = \text{(interval) constant}
\end{align*}
and the numerator
\begin{align*}
f(x) = f(x-y) + (y)f'(x) + \frac{1}{2}(y)^2f''(x) + \ldots \frac{1}{k!}(y)^{k}f^{k}(\eta),
\end{align*}

around $y = 0$. Here $\eta$ belongs to an intermediate point between $x$ and $x-y$, and we can enclose $f^{k}(\eta)$ in the whole interval $f^{k}([x-\varepsilon_1,x+\varepsilon_1]) \subset f^{k}([-\pi,\pi])$. Finally, we can factor out $(y)$ and divide both in the numerator and the denominator, getting

\begin{align*}
\frac{1}{\pi}\int_{|y| < \varepsilon_1} \frac{f'(x) + \frac{1}{2}(y)f''(x) + \ldots \frac{1}{k!}(y)^{k-1}f^{k}(\eta)}{1+cy^2}dy,
\end{align*}

which we could either bound or integrate explicitly since it is a smooth (interval) function and $f(x)$ is explicit. 

For $H^{F}f(x)$ we will do the same, expanding the cotangent function to avoid division by $\infty$:

\begin{align*}
\frac{1}{2}\cot\left(\frac{y}{2}\right) = -\frac{1}{4}(x-\pi) + c(\varepsilon_2)(x-\pi)^3, \quad c(\varepsilon_2) = \text{(interval) constant}
\end{align*}

we obtain

\begin{align*}
\frac{1}{\pi}\int_{|y-\pi| < \varepsilon_2} (f(x)-f(x-y))\left(-\frac{1}{4}(y-\pi) + c(\varepsilon_2)(y-\pi)^3\right)dy,
\end{align*}

which is smooth and therefore we can also integrate it as in the previous subsection.

The choice of $\varepsilon_i$ is determined by the balance between accuracy and computation time. Most of the times, the $\varepsilon_i$ will be taken very small and $H^{N}f(x)$ and $H^{F}f(x)$ will be regarded as error terms.

In the case where the integration domain is unbounded (for simplicity we may assume it is $\mathbb{R}$ and the integrand decays fast enough), one can do two workarounds:

\begin{itemize}
 \item Perform a change of variables that maps $\mathbb{R}$ onto a bounded domain, such as $x = 2\tan\left(\frac{y}{2}\right)$. This change of variables is useful because the problem is mapped onto $[-\pi,\pi]$ and one can work with Fourier series there. The integral in the new coordinates becomes

\begin{align*}
 \int_{-\infty}^{\infty} f(x)dx = \int_{-\pi}^{\pi}f\left(2\tan\left(\frac{y}{2}\right)\right) \sec\left(\frac{y}{2}\right)^{2}dy,
\end{align*}

which depending on $f$ may potentially be singular, in which case we would apply the techniques outlined in the beginning of this subsection.

\item Choose a large enough number $M$ and treat the contribution to the integral from $|x| > M$ as an error. Thus:

\begin{align*}
 I = \int_{-\infty}^{\infty} f(x)dx = \int_{|x| \leq M} f(x) dx + \int_{|x|>M} f(x) dx = I_1 + I_2.
\end{align*}

The term $I_1$ will be integrated normally. For the term $I_2$, assuming $|f(x)| \leq \frac{C}{|x|^{k}}$ we easily obtain $|I_2| \leq \frac{2C}{k-1} \frac{1}{M^{k-1}}$. Making $M$ large this term will go to zero.

\end{itemize}

\section{The Muskat Problem}\label{sectionmuskat}

The first problem that we will present is the so-called Muskat problem. This problem models the evolution of the interface between two different incompressible fluids with the same viscosity in a two-dimensional porous medium, and is used in the context of oil wells \cite{Muskat:porous-media}.

The setup consists of two incompressible fluids with different densities, $\rho^{1}$ and $\rho^{2}$, and the same viscosity, evolving in a porous medium with permeability $\kappa(x)$. The velocity is given by Darcy's law:

\begin{equation}\label{IIIdarcy}
\mu\frac{v}{\kappa}=-\nabla p-g\left(\begin{array}{cc}0 \\ \rho\end{array}\right),
\end{equation}
where $\mu$ is the viscosity, $p$ is the pressure and $g$ is the acceleration due to gravity, and the incompressibility condition
\begin{equation}\label{IIIincom}
\nabla\cdot v=0.
\end{equation}
We take $\mu=g=1$. The fluids also satisfy the conservation of mass equation
\begin{equation}\label{IIIconser}
\pat\rho+v\cdot\nabla\rho=0.
\end{equation}
We will work in different settings: the flat at infinity and horizontally periodic cases $(\Omega = \RR^{2}$ and $\Omega = \TT \times \RR$ respectively) or the confined case ($\Omega = \RR \times \left(-\frac{\pi}{2},\frac{\pi}{2}\right)$).   We denote by $\Omega^{1}$ the region occupied by the fluid with density $\rho^{1}$ (the ``top'' fluid) and by $\Omega^{2}$ the region occupied by the fluid with density $\rho^{2} \neq \rho^{1}$ (the ``bottom'' fluid).  All quantities with superindex $1$ (resp. $2$) will refer to $\Omega^{1}$ (resp. $\Omega^{2}$). The interface between the two fluids at any time $t$ is a planar curve denoted by $z(\cdot,t)$.

In the case $\Omega = \RR^{2}$, one can rewrite the system \eqref{IIIdarcy}--\eqref{IIIconser} in terms of the curve $z=(z^1,z^2)$, obtaining
\begin{align}\label{muskatinterface}
\partial_{t} z(\al,t) = \frac{\rho^{2} - \rho^{1}}{2\pi} P.V. \int_\mathbb{R} \frac{z^1(\al,t) - z^1(\be,t)}{|z(\al,t) - z(\be,t)|^{2}}(\partial_{\al}z(\al,t) - \partial_{\be}z(\be,t)) d\be.
\end{align}
In the horizontally periodic case ($\Omega = \TT \times \RR$) with $z(\alx+2\pi,t)=z(\alx,t)+(2\pi,0)$,  the evolution of the curve is given by the formula
\begin{align}\label{muskatinterfaceperiodic}
\partial_{t} z(\al,t) = \frac{\rho^{2} - \rho^{1}}{4\pi} \int_\mathbb{T} \frac{\sin(z^1(\al,t) - z^1(\be,t))(\partial_{\al}z(\al,t) - \partial_{\be}z(\be,t))}{\cosh(z^{2}(\al,t) - z^{2}(\be,t)) - \cos(z^{1}(\al,t) - z^{1}(\be,t))} d\be.
\end{align}

Finally, in the confined case:

\begin{multline*}
\pat z(\alpha,t)=\frac{\rho^{2}-\rho^{1}}{4\pi}\int_\RR\frac{(\paa z(\alpha,t)-\paa z(\alpha-\beta,t))\sinh(z_1(\alpha,t)-z_1(\alpha-\beta,t))}{\cosh(z_1(\alpha,t)-z_1(\alpha-\beta,t))-\cos(z_2(\alpha,t)-z_2(\alpha-\beta,t))}d\beta\\
+\frac{\rho^{2}-\rho^{1}}{4\pi}\int_\RR\frac{(\paa z_1(\alpha,t)-\paa z_1(\alpha-\beta,t),\paa z_2(\alpha,t)+\paa z_2(\alpha-\beta,t))\sinh(z_1(\alpha,t)-z_1(\alpha-\beta,t))}{\cosh(z_1(\alpha,t)-z_1(\alpha-\beta,t))+\cos(z_2(\alpha,t)+z_2(\alpha-\beta,t))}d\beta.
\end{multline*}

We define the Rayleigh-Taylor condition
$$
RT(\alpha,t)=-(\nabla p^2(z(\alpha,t))-\nabla p^1(z(\alpha,t)))\cdot\partial_\alpha^\bot z(\alpha,t),
$$

which can be written in terms of the interface as

\begin{align*}
 RT(\al,t) = (\rho^{-} - \rho^{+})\partial_{\al} z^{1}(\al,t).
\end{align*}

Linearizing around the steady state $(\alpha,0)$, the evolution equation of a small perturbation $(0,\varepsilon f_{L}(\al,t))$ satisfies at the linear level:

\begin{align}\label{linearf}
\pa_t f_{L}(\al,t) = -RT^{L}(\al,t) \Lambda f_{L}(\al,t),
\end{align}

where $RT^{L}(\al,t)$ is the linearized version of the Rayleigh-Taylor condition

$$
RT^{L}(\al,t)=g(\rho^2-\rho^1)
$$

 and $\Lambda = (-\Delta)^{\frac12}$. Thus, the equation \eqref{linearf} is parabolic if $RT^{L}(\al,t) > 0$. At the nonlinear level, similar estimates of the form

\begin{align*}
\pa_{t} \pa_{\al}^{k}z(\al,t) = -RT(\al,t) \Lambda \pa_{\al}^{k}z(\al,t) + \text{ lower order terms }
\end{align*}

can be derived for large enough $k$. It is therefore easy to see that the sign of $RT(\al,t)$ is crucial, since it will govern the stability of the equation (stable for positive $RT$, negative otherwise). For a fixed $t$, if $RT(\alpha,t)>0 , \;\forall\alpha\in\RR$ we will say that the curve is in the Rayleigh-Taylor stable regime and if $RT(\alpha,t)<0$ for some $\alpha$, we will say that the curve is in the Rayleigh-Taylor unstable regime. In other words, we have the correspondence:

\begin{align*}
z(\al,t) \text{ can be parametrized as a graph} \Leftrightarrow z(\al,t) \text{ is in the R-T stable regime}
\end{align*}

We will also say that if the curve $z(\al,t)$ changes from graph to non-graph or viceversa, it undergoes a stability shift.

\subsection{Brief history of the problem}

The Muskat problem  has been studied in many works. A proof of local existence of classical solutions in the Rayleigh-Taylor stable regime and a maximum principle for $\| \pa_xf(\cdot,t)\|_{L^\infty}$ can be found in \cite{Cordoba-Gancedo:maximum-principle-muskat}. See also \cite{Ambrose:well-posedness-hele-shaw}. Ill-posedness in the unstable regime appears in \cite{Cordoba-Gancedo:contour-dynamics-3d-porous-medium}.

Moreover, the authors in \cite{Cordoba-Gancedo:contour-dynamics-3d-porous-medium} showed that if $\|\partial_{x} f_0\|_{L^\infty}<1$, then $\|\partial_{x} f(\cdot,t)\|_{L^\infty}\le \|\partial_{x} f_0\|_{L^\infty}$ for all $t>0$. Further work has shown instant analyticity and existence of finite time turning \cite{Castro-Cordoba-Fefferman-Gancedo-LopezFernandez:rayleigh-taylor-breakdown}: in other words, the curve ceases to be a graph in finite time and the Rayleigh-Taylor condition changes sign to negative somewhere along the curve. The gap between these two results (i.e., the question whether the constant 1 is sharp or not for guaranteeing global existence) is still an open question, and there is numerical evidence of data with $\|\pa_x f_{0}\|_{L^{\infty}} = 50$ which turns over \cite{Cordoba-GomezSerrano-Zlatos:stability-shifting-muskat}.

 Given the parabolic character of the equation, it is natural to expect global existence, at least for small initial data. The first proof for small initial data was carried out in \cite{Siegel-Caflisch-Howison:global-existence-muskat} in the case where the fluids have different viscosities and the same densities (see \cite{Cordoba-Gancedo:contour-dynamics-3d-porous-medium} for the setting of the present paper --- different densities and the same viscosities --- and also \cite{Cheng-GraneroBelinchon-Shkoller:well-posedness-h2-muskat} for the general case).  The regularity requirement of the initial data has been subsequently lowered down: see \cite{Constantin-Cordoba-Gancedo-Strain:global-existence-muskat,Constantin-Cordoba-Gancedo-RodriguezPiazza-Strain:muskat-global-2d-3d,Constantin-Gancedo-Shvydkoy-Vicol:global-regularity-muskat-finite-slope,Matioc:local-existence-muskat-hs,Beck-Sosoe-Wong:duchon-robert-solutions-muskat,Cordoba-Lazar:global-wellposedness-muskat-H32,Cameron:global-wellposedness-muskat-slope-less-1,Deng-Lei-Lin:2d-muskat-monotone-data} for recent developments of global existence in different spaces.  A blow-up criterion was found in \cite{Constantin-Gancedo-Shvydkoy-Vicol:global-regularity-muskat-finite-slope}. For large time estimates see \cite{Patel-Strain:large-time-estimates-muskat,Gancedo-GarciaJuarez-Patel-Strain:muskat-viscosity-global}. In the case where surface tension is taken into account, global existence was shown in \cite{Escher-Matioc:parabolicity-muskat,Friedman-Tao:nonlinear-stability-muskat-capillary-pressure}.

Contrary to the previous intuition, there is also formation of singularities: there are initial data which start smooth, but once the Rayleigh-Taylor condition is not satisfied (i.e. they have ceased to be a graph), the smoothness of the curve may break down in finite time \cite{Castro-Cordoba-Fefferman-Gancedo:breakdown-muskat}. Another possibility is the appearance of a finite time self-intersection of the free boundary, either at a point (``splash'' singularity) or along an arc (``splat'' singularity). In the stable density jump case these cannot happen \cite{Gancedo-Strain:absence-splash-muskat-SQG,Cordoba-Gancedo:absence-squirt-singularities-muskat}. However, in the one phase case there are finite time splash singularities \cite{Castro-Cordoba-Fefferman-Gancedo:splash-singularities-muskat} but no splat singularities are possible \cite{Cordoba-PernasCastano:nonsplat-one-phase-muskat}.

 More general models, which take into account finite depth or non-constant permeability, and which also exhibit (single) stability shift were studied in \cite{Berselli-Cordoba-GraneroBelinchon:local-solvability-singularities-muskat,GomezSerrano-GraneroBelinchon:turning-muskat-computer-assisted,Cordoba-GraneroBelinchon-Orive:confined-muskat,GraneroBelinchon:global-existence-confined-muskat,PernasCastano:local-existence-inhomogeneous-muskat,Cordoba-PernasCastano:splash-splat-one-phase-muskat}.

For a bigger and more comprehensive list of references we refer the reader to the two surveys \cite{Castro-Cordoba-Gancedo:recent-results-muskat} and \cite{Gancedo:survey-muskat}.

\subsection{Theorems}

We now present a few theorems which can be proved by means of computer-assisted estimates, following the abstract idea from Example \ref{example1}. The intuition behind these results is driven by numerical simulations.

The first result compares the confined and flat at infinity regimes \cite{GomezSerrano-GraneroBelinchon:turning-muskat-computer-assisted}:

\begin{theorem}
There exists a family of analytic curves $z^{0}(\al) = (z_1^{0}(\al),z_2^{0}(\al))$, flat at infinity, for which there exists a finite time $T$ such that the solution to the confined Muskat problem develops a stability shift before $t=T$ and the non confined does not.

Specifically, we prove that the solution with initial data $z^{0}(\al)$ can be parametrized as a graph for all $t < T$ in the flat at infinity setting and can not be in the confined one.

\label{ThmConfTurnsNoConfNoTurns}
\end{theorem}

We now outline the reduction to a finite set of conditions.

\begin{reduction}{\ref{ThmConfTurnsNoConfNoTurns}}

We want to show that $RT(0,t) \sim At + O(t^2)$ for small, positive $t$, where $A < 0$ in the confined case and $A > 0$ in the flat at infinity, and this will ensure the Theorem. After some calculations one obtains:

\begin{multline*}
A_{confined}=2 \partial_{\alpha} z_2(0)\int_0^\infty \partial_\alpha z_1(\eta)\sinh(z_1(\eta))\sin(z_2(\eta))\bigg{(}\frac{1}{(\cosh(z_1(\eta))-\cos(z_2(\eta)))^2}\\
+\frac{1}{(\cosh(z_1(\eta))+\cos(z_2(\eta)))^2}\bigg{)}d\eta.
\end{multline*}
With the same approach, for the unconfined case the expression is
$$
A_{flat}=8 \partial_{\alpha} z_2(0)\int_0^\infty \frac{\partial_\alpha z_1(\eta)z_1(\eta)z_2(\eta)}{(z_1(\eta))^2+(z_2(\eta))^2)^2}d\eta.
$$

Thus, the theorem will be proved if we manage to validate the following open conditions:
\begin{align}
\label{signconfnoconf}
A_{confined}<0,\quad A_{flat}>0.
\end{align}

\end{reduction}

\begin{figure}[ht!]\centering
\includegraphics[scale=0.35]{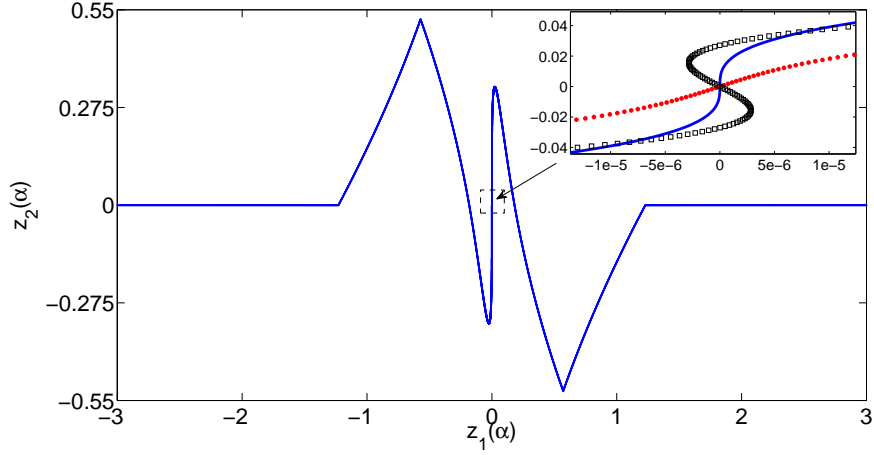}
\caption{The curve in Theorem \ref{ThmConfTurnsNoConfNoTurns}. Inset: Close caption around zero, solid: initial condition, dotted: normal component of the velocity for the flat at infinity case, squared: normal component of the velocity for the confined case.}
\label{figteo1}
\end{figure}

We can rigorously validate them for the following data (see Figure \ref{figteo1}):

\begin{align*}
z_1(\al) & = \al - \sin(\al)e^{-B\al^{2}}, \quad B = 10^{-4} \\
z_2(\al) & = \left\{
\begin{array}{lr}
\displaystyle \frac{\sin(3\al)}{3} & \displaystyle \text{ if } 0 \leq \al \leq \frac{\pi}{3} \vspace{0.2cm} \\
\displaystyle - \al + \frac{\pi}{3} & \displaystyle \text{ if } \frac{\pi}{3} \leq \al \leq \frac{\pi}{2} \vspace{0.2cm} \\
\displaystyle \al - \frac{2\pi}{3} & \displaystyle \text{ if } \frac{\pi}{2} \leq \al \leq \frac{2\pi}{3}\vspace{0.2cm} \\
\displaystyle 0 & \displaystyle \text{ if } \frac{2\pi}{3} \leq \al, \\
\end{array}
\right.
\end{align*}
where $z_2$ is extended such that it is an odd function. In Figure \ref{figteo1} (inset), we plot the normal velocity around the vertical tangent for the two scenarios (confined and flat at infinity), both scaled by a factor $1/100$. We can observe that the velocity denoted by squares, which corresponds to the confined case, will make the curve develop a turning singularity, where the dotted one (non-confined case) will force the curve to stay in the stable regime.

In fact, one can find more striking behaviours if the expansion of $RT(\al,t)$ is done at higher order (at the price of higher complexity of the calculations). The following Theorem was proved in \cite{Cordoba-GomezSerrano-Zlatos:stability-shifting-muskat-II}:

\begin{theorem}\label{theoremshifting}
There exist $T>\gamma>0$  and a spatially analytic solution $z$ to \eqref{muskatinterfaceperiodic} on the time interval $[-T,T]$  such that $z(\cdot,t)$ is a graph of a smooth function of $\al$ when $|t|\in[T-\gamma,T]$ (i.e., $z$ is in the stable regime near $t=\pm T$)  but $z(\cdot,t)$ is not a graph of a function of $x$ when $|t|\le \gamma$ (i.e., $z$ is in the unstable regime near $t=0$).

In other words, there exists solutions of \eqref{muskatinterfaceperiodic} that make the transition stable $\rightarrow$ unstable $\rightarrow$ stable.
\end{theorem}

The intuition behind this result comes from the numerical experiments which were started in \cite{Cordoba-GomezSerrano-Zlatos:stability-shifting-muskat}, where it was proved that there were solutions that were exhibiting the unstable $\rightarrow$ stable $\rightarrow$ unstable transition. These suggested  existence of curves which are (barely) in the unstable regime, and such that the evolution both forward and backwards in time transports them into the stable regime. (We note that neither the velocity nor any other quantity was observed to become degenerate  in these experiments). We remark that this behaviour is purely nonlinear and thus nonlinear effects may dominate the linear ones under certain conditions. 

\begin{reduction}{\ref{theoremshifting}}
  Let $\ep\ge 0$ and consider the initial family of curves $z_{\ep}(\al,0)=(z^{1}_{\ep}(\al,0),z^{2}_{\ep}(\al,0))$, with
\begin{align*}
 z^{1}_{\ep}(\al,0) & = \al - \sin(\al) - \ep \sin(\al), \\
z^{2}_{\ep}(\al,0) & = A(\ep) \sin(2\al).
\end{align*}

The goal is to show that this family of solutions satisfies $RT(0,t) \sim -\varepsilon + Ct^2 + O(t^3)$. See Figure \ref{initialcondition}.

\begin{figure}[ht]
\centering
\includegraphics[width=\textwidth]{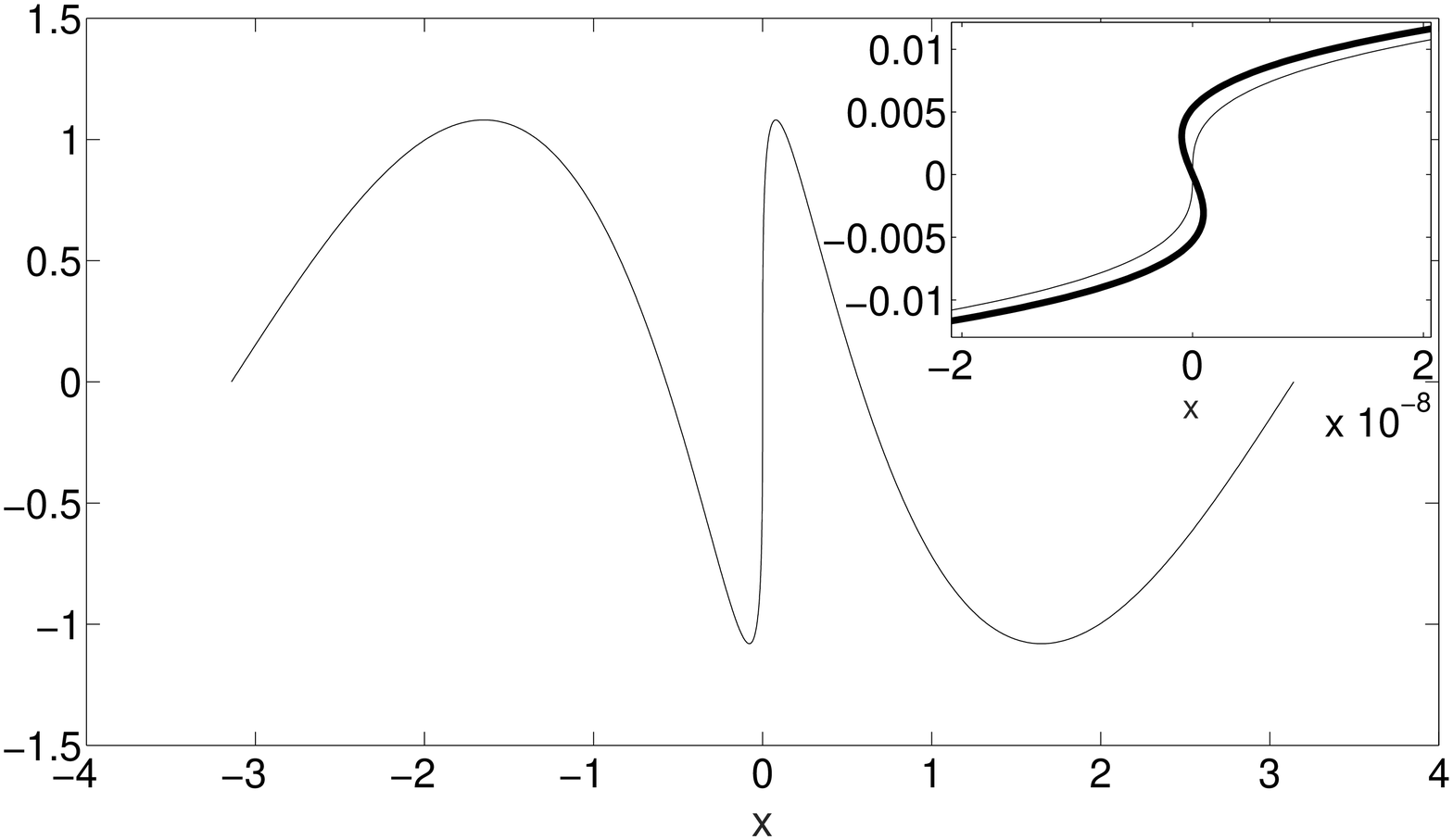}
\caption{$z_{\ep}(\al,0)$ from Theorem \ref{theoremshifting} with $A(\ep) = 1.08050$. Inset: closeup around $x = 0$. Thick curve: $\ep = 10^{-6}$, thin curve: $\ep = 0$. We remark that both curves are indistinguishable at the larger scale.}
\label{initialcondition}
\end{figure}

This is done in two steps. The first one is to choose $A(\ep)$ accordingly: one can prove that for any $\ep\in[0,10^{-6}]$, there exists $A(\ep)\in(1.08050, 1.08055)$ such that if $z_\ep$ solves \eqref{muskatinterfaceperiodic} with  initial data $z_{\ep}(\al,0)$, then
\begin{align*}
 \pa_{t} RT(0,0) = 0.
\end{align*}

The second one is to show that there exist $T > 0,C \geq 1$, independent of $\ep$, such that for any $\ep\in[0,10^{-6}]$ and $A(\ep)$ chosen before, there is a unique analytic solution $z_{\ep}$ of \eqref{muskatinterfaceperiodic} on the time interval $(-T,T)$ with initial data $z_{\ep}(\al,0)$, and it satisfies
\begin{equation} \label{2.1}
\pa_{tt} RT(0,0)\ge 30.
\end{equation}

The first step is accomplished calculating $\pa_{t} RT(0,0)$ for $A(\ep) = 1.08050$ and for $A(\ep) = 1.08055$ and checking that one is negative and the other is positive, for each $\ep\in[0,10^{-6}]$. The second step follows by taking $A(\ep) = [1.08050,1.08055]$ (the full interval, since we do not know $A(\ep)$ explicitly) and propagating this interval in the relevant computations. The drawback of this method is that $\pa_{tt} RT(0,0)$ consists of tens of terms of the type

\begin{multline*}
 B_{11}(\al) = -\int_{\T} \int_{\T} \frac{\sin(z^1(\al) - z^1(\al-y)))(z^{1}_{\al}(\al) - z^{1}_{\al}(\al-y))^{2}}{\cosh(z^2(\al) - z^2(\al-y)) - \cos(z^1(\al) - z^1(\al-y))} \\
 \times \left(\frac{\sin(z^1(\al) - z^1(\al-z))(z^{1}_{\al}(\al) - z^{1}_{\al}(\al-z))}{\cosh(z^2(\al) - z^2(\al-z)) - \cos(z^1(\al) - z^1(\al-z))}\right.
\\
\left. - 
\frac{\sin(z^1(\al-y) - z^1(\al-y-z))(z^{1}_{\al}(\al-y) - z^{1}_{\al}(\al-y-z))}{\cosh(z^2(\al-y) - z^2(\al-y-z)) - \cos(z^1(\al-y) - z^1(\al-y-z))}
\right)
dy dz 
\end{multline*}

which are 2-dimensional integrals which contain a singularity. In order to overcome this issue, natural extensions (to 2D) of the schemes outlined in subsection \ref{subsection:singularint} in the 1D case need to be done.

\end{reduction}

\begin{corollary}
 Using the same techniques as in \cite{Cordoba-GomezSerrano-Zlatos:stability-shifting-muskat}, we can construct solutions that change stability 4 times according to the transition unstable $\rightarrow$ stable $\rightarrow$ unstable $\rightarrow$ stable $\rightarrow$ unstable.
\end{corollary}

\begin{figure}[ht!]\centering
\includegraphics[scale=0.35]{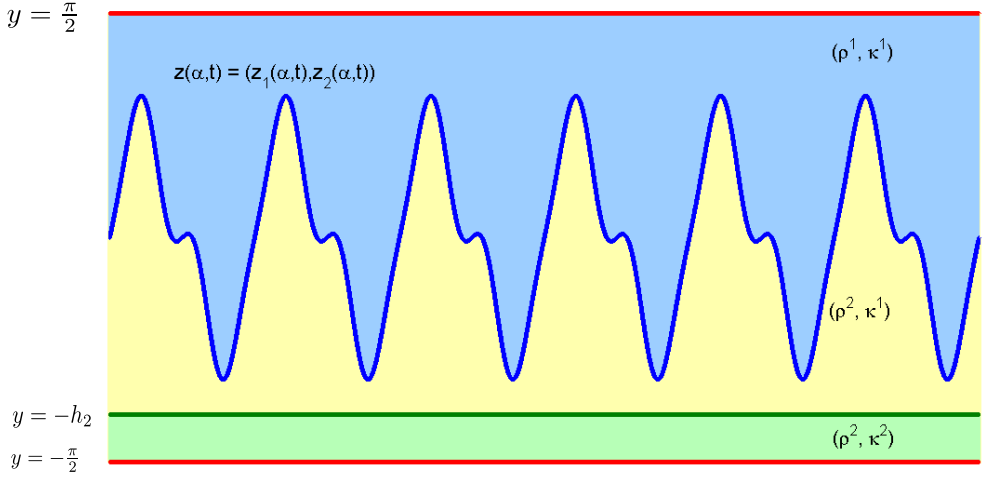}
\caption{Situation of the different fluids and permeabilities.We have 3 regions with parameters $(\rho_1, \kappa_1), (\rho_2,\kappa_1),(\rho_2,\kappa_2)$ separated by two boundaries $z(\al,t)$ and $(\al,-h_2)$ in a confined medium. The first one is a free boudnary and the second is fixed.}
\label{FigInhomogeneous}
\end{figure}

The last theorem we present applies to a model where we also take into account a permeability jump. This problem is important in the context of geothermal reservoirs \cite{Cerminara-Fasano:dynamics-geothermal-reservoir}, where it could represent different types of rock layers (impermeable, permeable), below which a heat source (magma) is located. See Figure \ref{FigInhomogeneous} for a depiction of the setting.

In that case, the evolution equation for the interface $z(\al,t)$ is more complicated and given by:

\begin{align}
\pat z(\alpha)=&\bar{\rho}\text{P.V.}\int_\RR\frac{(\paa z(\alpha)-\paa z(\beta))\sinh(z_1(\alpha)-z_1(\beta))}{\cosh(z_1(\alpha)-z_1(\beta))-\cos(z_2(\alpha)-z_2(\beta))}d\beta\nonumber\\
&+\bar{\rho}\text{P.V.}\int_\RR\frac{(\paa z_1(\alpha)-\paa z_1(\beta),\paa z_2(\alpha)+\paa z_2(\beta))\sinh(z_1(\alpha)-z_1(\beta))}{\cosh(z_1(\alpha)-z_1(\beta))+\cos(z_2(\alpha)+z_2(\beta))}d\beta\nonumber\\
&+\frac{1}{4\pi}\text{P.V.}\int_\RR\varpi_2(\beta)BS(z_1(\alpha),z_2(\alpha),\beta,-h_2)d\beta\nonumber\\
&+\frac{\paa z(\alpha)}{4\pi}\text{P.V.}\int_\RR\varpi_2(\beta)\frac{\sin(z_2(\alpha)+h_2)}{\cosh(z_1(\alpha)-\beta)-\cos(z_2(\alpha)+h_2)}d\beta\nonumber\\
&+\frac{\paa z(\alpha)}{4\pi}\text{P.V.}\int_\RR\varpi_2(\beta)\frac{\sin(z_2(\alpha)-h_2)}{\cosh(z_1(\alpha)-\beta)+\cos(z_2(\alpha)-h_2)}d\beta\label{IIIeqv2}.
\end{align}

with

\begin{align}
\varpi_2(\al)=&-2\mathcal{K}BR(\varpi_1,z)h(\al)\cdot(1,0)+\frac{2\mathcal{K}^2}{2\pi}BR(\varpi_1,z)h(\al)\cdot(1,0)*G_{h_2,\mathcal{K}}\nonumber\\
=&\, 2\mathcal{K}\bar{\rho}\left[\text{P.V.}\int_\RR\paa z_2(\beta)\frac{\sin(h_2+z_2(\beta))}{\cosh(\alpha-z_1(\beta))-\cos(h_2+z_2(\beta))}d\beta\nonumber\right.\\
&-\text{P.V.}\int_\RR\paa z_2(\beta)\frac{\sin(-h_2+z_2(\beta))}{\cosh(\alpha-z_1(\beta))+\cos(-h_2+z_2(\beta))}d\beta\nonumber\\
&-\frac{\mathcal{K}}{2\pi}G_{h_2,\mathcal{K}}*\text{P.V.}\int_\RR\frac{\paa z_2(\beta)\sin(h_2+z_2(\beta))}{\cosh(\alpha-z_1(\beta))-\cos(h_2+z_2(\beta))}d\beta\nonumber\\
&+\left.\frac{\mathcal{K}}{2\pi}G_{h_2,\mathcal{K}}*\text{P.V.}\int_\RR\frac{\paa z_2(\beta)\sin(-h_2+z_2(\beta))}{\cosh(\alpha-z_1(\beta))+\cos(-h_2+z_2(\beta))}d\beta\right], \label{IIIw2defb}
\end{align}

and

$$
G_{h_2,\mathcal{K}}(\xi)=\int_\RR\frac{\cos(y\xi)\sinh(2h_2 y)}{\sinh(\pi y)+\mathcal{K}\sinh(2h_2 y)}dy, \quad \mathcal{K}=\frac{\kappa^1-\kappa^2}{\kappa^1+\kappa^2}, \quad \bar{\rho}=\frac{\kappa^1(\rho^2-\rho^1)}{4\pi}.$$

\begin{figure}[ht!]\centering
\includegraphics[trim=2cm 0mm 0mm 0mm,scale=0.55]{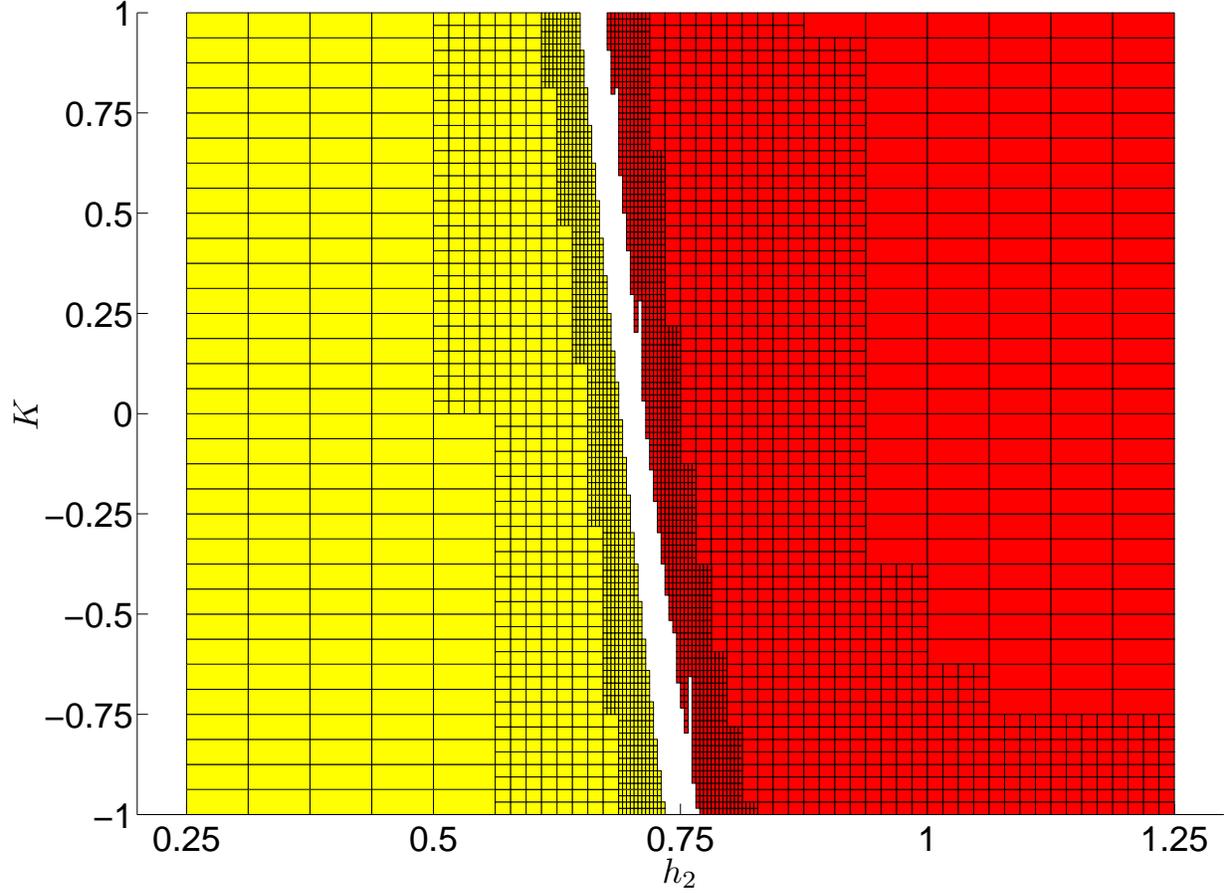}
\caption{Bifurcation diagram corresponding to the phenomenon of stability shift for the initial condition given by the family of curves \eqref{initialConditionBif}. Yellow (lighter color): no shift, red (darker color): shift. The boundary separating the two colors is smooth and can be parametrized as $(h_2,\mathcal{K}(h_2))$.}
\label{FigBifurcacion}
\end{figure}

The goal is to illustrate the different behaviours that may arise by looking at the short-time evolution of a family of initial data, depending on the height of the permeability jump and the magnitude of the permeabilities. This is shown in the bifurcation diagram in Figure \ref{FigBifurcacion}, where we plot whether for short time, the curve will shift stability or not.

\begin{theorem}
\label{ThmConfinedInhomogeneous}
There exists a family of analytic initial data $z(\al,h_2) = (z_1(\al,h_2),z_2(\al,h_2))$, depending on the height at which the permeability jump is located, such that the corresponding solution to the confined, inhomogeneous Muskat \eqref{IIIeqv2}-\eqref{IIIw2defb}:
\begin{enumerate}[(a)]
\item
\begin{enumerate}[1.]
\item For all $0.25 < h_2 < h_2^{ntu} = 0.648$, the curve will not shift independently of $\mathcal{K}$.
\item For all $ 0.676 < h_2 < 0.686$, the permeabilities help the shift.
\item For all $ 0.715 < h_2 <  0.738$, the permeabilities prevent the shift.
\item For all $ 0.77 = h_2^{tu} < h_2 < 1.25$, the curve will shift independently of $\mathcal{K}$.
\end{enumerate}
\item There exists a $C^1$ curve $(h_2,\mathcal{K}(h_2))$, located in $[0.648,0.77] \times (-1,1)$, such that for every $h_2$ for which the curve is defined, for every $\mathcal{K}<\mathcal{K}(h_2)$ the curve does not turn and for every $\mathcal{K}>\mathcal{K}(h_2)$ the curve turns.
\end{enumerate}
\end{theorem}

\begin{reduction}{\ref{ThmConfinedInhomogeneous}(a)}

We proceed now to calculate $\pa_{t} RT(0,0)$. Then, the appropriate expression is
$$
\pa_{t} RT(0,0)=C(I_1+I_2),
$$

for some $C > 0$, where
$$
I_1=2\paa z_2(0)\int_0^\infty \frac{\paa z_1(\beta)\sinh(z_1(\beta))\sin(z_2(\beta))}{\left(\cosh(z_1(\beta))-\cos(z_2(\beta))\right)^2}+\frac{\paa z_1(\beta)\sinh(z_1(\beta))\sin(z_2(\beta))}{\left(\cosh(z_1(\beta))+\cos(z_2(\beta))\right)^2}d\beta,
$$
and

\begin{multline*}
I_2=4\paa z_2(0)\mathcal{K}\int_0^\infty\int_0^\infty
\frac{\paa z_2(\gamma)\cos(z_1(\gamma)y)}{(\sinh(\pi y)+\mathcal{K}\sinh(2h_2 y))\cosh\left(y\frac{\pi}{2}\right)}\\
\times\left(2y\cosh\left(\frac{y\pi}{2}- y h_2\right)\cosh\left(\frac{y \pi}{2}\right)-\frac{2\sinh\left(y h_2\right)}{\tan(h_2)}\right)\\
\times \cosh\left(y z_2(\gamma)\right)\cosh\left(y\left(\frac{\pi}{2}-h_2\right)\right) d\gamma dy.
\end{multline*}

\end{reduction}

The initial condition family we used for the bifurcation diagram was
\begin{align}
z_1(\al) & = \al - \sin(\al)e^{-B \al^{2}}, \quad B = 10^{-4} \nonumber \\
z_2(\al) & = h_2\frac{3}{\pi}\left(\frac{\sin(3\al)}{3}- \frac{\sin(\al)}{2.5}\left(e^{-(\al+2)^2}+e^{-(\al-2)^2}\right)\right)1_{\{|\al| \leq \pi\}}. \label{initialConditionBif}
\end{align}

We computed the bifurcation diagram depicted in Figure \ref{FigBifurcacion}. We could give an answer regarding the question of short time stability shifting to $97.14\%$ of the parameter space. $53.23\%$ of the space turned (red) and $43.91\%$ did not turn (yellow). The remaining $2.86\%$ is painted in white.

We proceeded as follows: for each region in parameter space, we computed an enclosure of $\pa_t RT(0,0$) for all values in that region. If we could establuish a sign we painted the region of the corresponding color. If not, we subdivided the region and recomputed up to a certain maximum number of subdivisions. We remark that due to the actual answer being very close to zero or zero, the enclosures were not conclusive in some regions and more precision is required for those.

\begin{reduction}{\ref{ThmConfinedInhomogeneous}(b)}
 We want to invoke the Implicit Function Theorem. Thus, we have to check that 
$$
\frac{d}{d\mathcal{K}}\pat RT(0,0)\neq0\text{ for points }(h_2,\mathcal{K}) \text{ such that }\pat RT(0,0)=0.
$$
In particular, we have to check the previous condition in an open set containing the white region in Figure \ref{FigBifurcacion}. We compute
\begin{multline*}
DI_2 \equiv \frac{d}{d\mathcal{K}}\pat RT(0,0)=4\paa z_2(0)\int_0^\infty\int_0^\infty
\frac{\sinh(\pi y)\paa z_2(\gamma)\cos(z_1(\gamma)y)}{(\sinh(\pi y)+\mathcal{K}\sinh(2h_2 y))^2\cosh\left(y\frac{\pi}{2}\right)}\\
\times\left(2y\cosh\left(\frac{y\pi}{2}- y h_2\right)\cosh\left(\frac{y \pi}{2}\right)-\frac{2\sinh\left(y h_2\right)}{\tan(h_2)}\right)\\
\times \cosh\left(y z_2(\gamma)\right)\cosh\left(y\left(\frac{\pi}{2}-h_2\right)\right) d\gamma dy.
\end{multline*}

and show it is always non-zero.
\end{reduction}

\section{The surface quasi-geostrophic equation}\label{sectionsqg}

The Surface Quasi-Geostrophic equation (SQG) is the following active scalar equation
\begin{eqnarray}\label{equationsqg}
\left (\partial_t + u\cdot\nabla \right )\theta = 0
\end{eqnarray}
 where the relation between the incompressible velocity $u$ and $\theta$ is given by 
\begin{align*}
 u = \nabla^{\perp}\psi, \quad \theta=-(-\Delta)^{\frac{1}{2}}\psi.
\end{align*}
The scalar $\theta(x,t)$  represents the temperature and $\psi(x,t)$ is the stream function. The
 non-local operator $\Lambda^{\gamma} = (-\Delta)^{\frac{\gamma}{2}}$ is defined
through the Fourier transform by $\widehat{\Lambda^{\gamma}f}(\xi)=|\xi|^{\gamma}\hat{f}(\xi)$.  

This equation has applications to meteorology and oceanography, since it comes from models of atmospheric and ocean fluids \cite{Pedlosky:geophysical}
and is a special case of the more general 3D quasi-geostrophic
equation. There was a high scientific interest to understand
the behavior of the SQG equation, initially because it is a plausible model to
explain the formation of fronts of hot and cold air \cite{Pedlosky:geophysical}, and more recently \cite{Constantin-Majda-Tabak:formation-fronts-qg}
this system was proposed as a 2D model for the 3D vorticity
intensification and a geometric and
analytic analogy with the 3D incompressible Euler equations was shown.

One can also see a strong analogy with the 2D Euler equation in vorticity form:

\begin{eqnarray}\label{eulerequation}
\left (\partial_t + u\cdot\nabla \right )\omega = 0
\end{eqnarray}
\begin{align*}
 u = \nabla^{\perp}\psi, \quad \omega=-(-\Delta)^{1}\psi,
\end{align*}

the only difference being a stronger singular character of the velocity in the SQG case.

The problem of whether the SQG system presents finite time singularities or there is global existence is open for the smooth case. 

\subsection{Brief history of the problem}

Local existence has been proved in various functional settings \cite{Constantin-Majda-Tabak:formation-fronts-qg,Chae:qg-equation-triebel-lizorkin,Li:existence-theorems-2d-sqg-plane-waves,Wu:qg-equations-morrey-spaces,Wu:solutions-2d-qg-holder}. Starting from initial data with infinite energy, a gradient blowup may occur \cite{Castro-Cordoba:infinite-energy-sqg}. For finite energy initial data, solutions may start arbitrarily small and grow arbitrarily big in finite time \cite{Kiselev-Nazarov:simple-energy-pump-sqg}.

The numerical simulations in \cite{Constantin-Majda-Tabak:formation-fronts-qg} proposed a blowup scenario in the form of a closing hyperbolic saddle. This was ruled out in \cite{Cordoba:nonexistence-hyperbolic-blowup-qg,Cordoba-Fefferman:growth-solutions-qg-2d-euler}. More modern numerical simulations were able to resolve past the initially predicted singular time and found no singularities \cite{Constantin-Lai-Sharma-Tseng-Wu:new-numerics-sqg}. A new scenario was proposed in in \cite{Scott:scenario-singularity-quasigeostrophic}, starting from elliptical configurations, that develops filamentation and after a few cascades, blowup of $\nabla \theta$.

Global existence of weak solutions in $L^{2}$ was shown in  \cite{Resnick:phd-thesis-sqg-chicago}, and extended to the class of initial data belonging to $L^{p}$ with $p > 4/3$ in \cite{Marchand:existence-regularity-weak-solutions-sqg}. Non-uniqueness of weak solutions has been proved in \cite{Buckmaster-Shkoller-Vicol:nonuniqueness-sqg}. See also \cite{Nahmod-Pavlovic-Staffilani-Totz:global-invariant-measures-gsqg}.

Through a different motivation, \cite{Cordoba-Fefferman-Rodrigo:almost-sharp-fronts-sqg,Fefferman-Rodrigo:almost-sharp-fronts-sqg,Fefferman-Luli-Rodrigo:spine-sqg-almost-sharp-front} the existence of a special type of solutions that are known as ``almost sharp fronts'' was studied. These solutions can be thought of as a regularization of a front, with a small strip around the front in which the solution changes (reasonably) from one value of the front to the other. See \cite{Gravejat-Smets:travelling-waves-smooth-sqg} for a construction of traveling waves.

\subsection{Main Theorem}

The main Theorem is the following \cite{Castro-Cordoba-GomezSerrano:global-smooth-solutions-sqg}:

\begin{theorem}
\label{globalsqg}
 There is a nontrivial global smooth solution for the SQG equations that has finite energy, is compactly supported and is 3-fold.
\end{theorem}

\begin{figure}[ht!]\centering
\includegraphics[scale=0.35]{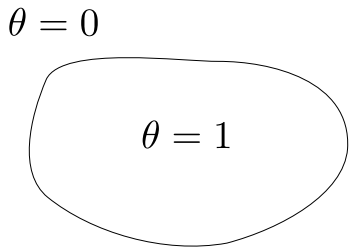}
\caption{The patch setting: the scalar $\theta$ is an indicator function of a domain $\Omega$. This setting is preserved in time: $\theta$ will be an indicator function of a moving domain $\Omega(t)$ for all $t$.}
\label{patch}
\end{figure}

It is well known that radial functions are stationary solutions of \eqref{equationsqg} due to the structure of the nonlinear term. The solutions that will be constructed are a smooth perturbation in a suitable direction of a specific radial function. The smooth profile we will perturb satisfies (in polar coordinates)
\begin{align*}
\theta(r)\equiv \left\{ \begin{array}{ccc} 1 & \text{for $0\leq r \leq 1-a$}\\ \text{smooth and decreasing} & \quad \text{for $1-a < r< 1$}  \\ 0 & \text{for $1\leq r <\infty$} \end{array}\right.,
\end{align*}
where $a$ is a small number. In addition the dynamics of these solutions consist of global rotating level sets with constant angular velocity. These level sets are a perturbation of the circle. The motivation comes from the so-called ``patch'' problem: namely when $\theta$ is a step function (see Figure \ref{patch}). In this setting, the uniformly rotating solutions are known as V-states.

In this setting, local existence of patch solutions has been obtained in \cite{Rodrigo:evolution-sharp-fronts-qg,Gancedo:existence-alpha-patch-sobolev} and uniqueness in \cite{Cordoba-Cordoba-Gancedo:uniqueness-sqg-patch}. There are two scenarios suggesting finite time singularities: the first one \cite{Cordoba-Fontelos-Mancho-Rodrigo:evidence-singularities-contour-dynamics}, starting from two patches, suggests an asymptotically self-similar collapse between the two patches, and at the same time a blowup of the curvature at the touching point; the second one \cite{Scott-Dritschel:self-similar-sqg} evolves a thin elliptical patch and indicates a self-similar filamentation cascade ending at a singularity with a blowup of the curvature.

The first computations of V-states for the 2D Euler equation are numerical \cite{Deem-Zabusky:vortex-waves-stationary} and since then there have been many works in different settings \cite{Wu-Overman-Zabusky:steady-state-Euler-2d,Elcrat-Fornberg-Miller:stability-vortices-cylinder,LuzzattoFegiz-Williamson:efficient-numerical-method-steady-uniform-vortices,Saffman-Szeto:equilibrium-shapes-equal-uniform-vortices}. However, the first proof dates to \cite{Burbea:motions-vortex-patches} proving their existence and later \cite{Hmidi-Mateu-Verdera:rotating-vortex-patch} their regularity. See also 
\cite{Hmidi-delaHoz-Mateu-Verdera:doubly-connected-vstates-euler,Hmidi-Mateu-Verdera:rotating-doubly-connected-vortices,Hmidi-Mateu:bifurcation-kirchhoff-ellipses,Hmidi-Mateu:degenerate-bifurcation-vstates-doubly-connected-euler,Hmidi-Renault:existence-small-loops-doubly-connected-euler,delaHoz-Hassainia-Hmidi-Mateu:vstates-disk-euler} for other studies in different directions (regularity of the boundary, different topologies, etc.).

For the SQG equation, in \cite{Castro-Cordoba-GomezSerrano:existence-regularity-vstates-gsqg,Castro-Cordoba-GomezSerrano:analytic-vstates-ellipses} existence and analyticity of the boundary was shown. In \cite{GomezSerrano:stationary-patches} nontrivial stationary solutions are constructed. The existence of doubly-connected V-states with analytic boundary was done in \cite{Renault:relative-equlibria-holes-sqg}. For more results concerning V-States of other active scalar equations see \cite{Hassainia-Hmidi:v-states-generalized-sqg,delaHoz-Hassainia-Hmidi:doubly-connected-vstates-gsqg,Hmidi-Mateu:existence-corotating-counter-rotating,Dritschel-Hmidi-Renault:imperfect-bifurcation-qgsw}.

We have also managed to prove a similar result as Theorem \ref{globalsqg} for other active scalar equations such as the 2D Euler equations \cite{Castro-Cordoba-GomezSerrano:uniformly-rotating-smooth-euler} (see also \cite{Garcia-Hmidi-Soler:non-uniform-vstates-euler}).

\begin{reduction}{\ref{globalsqg}}

We start writing down the scalar $\theta$ in terms of the level sets: $\theta(z(\al,\rho,t),t) = f(\rho)$ and write down the level sets in polar coordinates as $z(\al,\rho,t) = R(\lambda t)r(\al,\rho)(\cos(\al),\sin(\al))$ with $R$ being a rotation matrix. Then after a few algebraic manipulations one can show that a solution of 
the SQG equation that rotates with angular velocity $\lambda$ has to satisfy $F(r(\al,\rho),\lambda) = 0$, where $F(r(\al,\rho),\lambda)$ is 
an integrodifferential equation, and $r(\al,\rho)$ is the radial component of the solution. In this formulation, the coordinate $\rho$ is related to the level of the level sets, and $\alpha$ is an angular coordinate. The strategy is to apply an abstract Crandall-Rabinowitz theorem \cite{Crandall-Rabinowitz:bifurcation-simple-eigenvalues} bifurcating from $r = \rho$, which corresponds to 
a radial function (and therefore a solution for every $\lambda$). We bifurcate from a smooth ``annular'' profile which is 1 inside the disk of radius 0.95, 0 outside the disk of radius 1 and $C^4$ in between. 
The main steps to be shown are the following:
\begin{enumerate}
 \item[1.] $F$ is well defined and is $C^{1}$.

 \item[2.] Ker($\mathcal{F}$) is one-dimensional, where $\mathcal{F}$ is the linearized operator around $r = \rho$ at $\lambda = \lambda_{3}$, where $\lambda_{3}$ has to be determined.
 \item[3.]  $Y$/Range($\mathcal{F}$) is one-dimensional and $F_{r \lambda}(r,\lambda_{3})(r_3) \not \in$ Range($\mathcal{F}$), where Ker$(\mathcal{F}) = \langle r_3 \rangle$.
\end{enumerate}
The hardest step is step 2. We look for solutions with $3$-fold symmetry. To do so, we decompose the space and project onto the $3$-rd Fourier mode in $\alpha$. We try to find a function $B_{3}(\rho)$ in such a way that the kernel of $\mathcal{F}$ is generated by $\rho B_{3}(\rho) \cos(3\al)$. After rewriting the equations we end up having to solve an equation of the type
\begin{align*}
 \lambda_{3} B_{3}(\rho) = I(\rho)B_{3}(\rho) + \int T^{3}(\rho,\rho')B_{3}(\rho')d \rho',
\end{align*}

where both $I$ and $T^{3}$ are

\begin{align*}
I(\rho) & = -\frac{1}{2\pi} \int_{0}^{1} f_\rho(\rho')\left(\intpi \frac{\cos(x)}{\sqrt{1+\left(\frac{\rho}{\rho'}\right)^2-2\left(\frac{\rho}{\rho'}\right)\cos(x)}} dx\right) d\rho' \\
T^3(\rho,\rho') & = \frac{1}{2\pi} f_\rho(\rho') \frac{\rho'}{\rho}\int_{-\pi}^{\pi} \frac{\cos(mx)}{\sqrt{1+\left(\frac{\rho}{\rho'}\right)^2-2\left(\frac{\rho}{\rho'}\right)\cos(x)}}dx d\rho'.
\end{align*}

Note that these functions can also be written in terms of elliptic integrals. We regard the equation as an eigenvalue problem, having to find an eigenpair ($\lambda_{3}, B_{3}$). In fact, we will look for the smallest eigenvalue $\lambda_{3}$. The first drawback is that the integral operator is not symmetric: only close to symmetric so a priori it is not clear whether there is a (real) solution or not. Moreover, the appearance of the multiplication operator $I$ has the effect that the RHS is not compact, making this step more challenging. Nonetheless, we can prove the existence of $\lambda_{3}$ and $B_{3}$. 

The strategy is to consider this problem as a perturbation of a symmetric one. The hope is that if the antisymmetric part is small enough (compared to the gap between the eigenvalues), then there will be a real eigenpair. See Figure \ref{figspectrum} for a sketch of the situation: since there is only one eigenvalue inside the grey ball -- which is sufficiently small --, it has to be real (otherwise there would be two). We can recast it into explicit, quantitative conditions involving the gap between the first eigenvalues and the norm of the antisymmetric part.

\begin{figure}[ht]
\centering
\subfigure[Spectrum of the symmetric part of the RHS. $\lambda^{*}$ is the first eigenvalue, $c^*$ is the second. The rest of the spectrum is to the right of $c^*$.]
{
\includegraphics[width=0.65\textwidth]{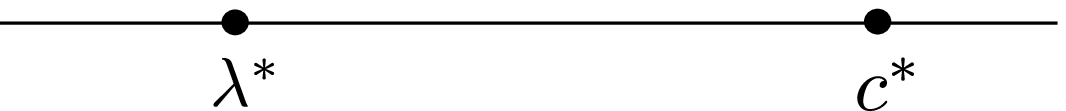}
\label{eigenvaluesym}
}
\subfigure[Spectrum of the RHS (perturbation of the symmetric part). There is an eigenvalue inside the grey ball and the rest of the spectrum lies inside the white ball or further right.]
{
\includegraphics[width=0.65\textwidth]{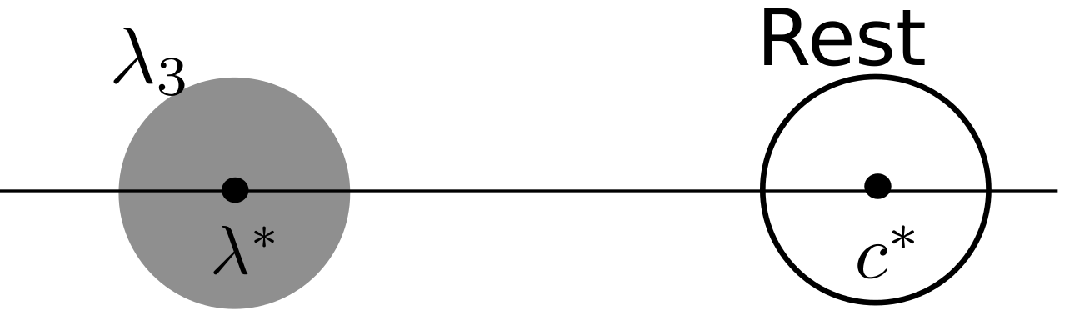}
\label{eigenvaluesymanti}
}
\caption{Sketch of the spectrum of the symmetric part of the RHS and the RHS.}
\label{figspectrum}
\end{figure}

We now focus on the symmetric part of the RHS. Getting a lower bound of the first eigenvalue is easy via Rayleigh-Ritz bounds. The hard part is to get a lower bound of the second eigenvalue of a symmetric operator.

We can get advantage of the compact part. If it were finite dimensional, then the problem reduces to finding an eigenvalue of a matrix. The crucial observation is that since the operator is compact, it will be well approximated by a finite rank operator modulo a small error, which can be made arbitrarily small increasing the dimension of the finite rank operator. The error can be written as an explicit (singular) integral which can be bounded as the ones in subsection \ref{subsection:singularint}. To deal with the singularity, we need to perform Taylor approximations of the elliptic functions with explicit error estimates. The philosophy can be summarized by the following: \textit{if we have an approximate guess of the eigenpair, then there is a true eigenpair nearby}. This way we can get tight, explicit bounds of the spectrum which lead to the proof of Step 2.

Step 1 is standard and technical, and step 3 follows from a similar analysis of the adjoint problem. 

\end{reduction}

\section*{Acknowledgements}
 J.G.-S. was partially supported by the grant MTM2014-59488-P (Spain), by the ICMAT-Severo Ochoa grant SEV-2015-0554, by the Simons Collaboration Grant 524109 and by the NSF-DMS 1763356 Grant. We would like to thank Diego C\'ordoba, Jordi-Llu\'is Figueras and Francisco Gancedo for helpful comments on previous versions of this manuscript.
 
This paper was developed out of a talk given at the XVIII Spanish-French School Jacques-Louis Lions about Numerical Simulation in Physics and Engineering, where I was awarded the 2018 Antonio Valle Prize from the Sociedad Espa\~nola de Matem\'atica Aplicada (SeMA). I would like to thank the SeMA and the organizers of the conference for such a great opportunity.

 \bibliographystyle{abbrv}
 \bibliography{references}

 \begin{tabular}{l}
 \textbf{Javier G\'omez-Serrano} \\
 {\small Department of Mathematics} \\
 {\small Princeton University}\\
 {\small 610 Fine Hall, Washington Rd,}\\
 {\small Princeton, NJ 08544, USA}\\
  {\small Email: jg27@math.princeton.edu}
 \end{tabular}

\end{document}